\documentclass[preprint, 12pt, draft]{elsarticle}

\usepackage[centertags]{amsmath}
\usepackage{multirow,newlfont,amssymb,amsfonts,amsthm ,booktabs, empheq, rotating, tablefootnote}

\DeclareMathOperator{\Span}{span}

\newdefinition{defn}{Definition}[section]

\theoremstyle{remark}

\newcommand{\p}{\partial}
\newcommand{\bb}{\begin{equation}}
\newcommand{\ee}{\end{equation}}
\newcommand{\ba}{\begin{array}}
\newcommand{\ea}{\end{array}}
\newcommand{\R}{\mathbb{R}}
\newcommand{\f}{\frac}
\usepackage[all]{xy}
\usepackage[usenames,dvipsnames]{color}
\usepackage{graphicx}
\usepackage{caption}
\usepackage{subcaption, longtable}

\makeatletter
\def\ps@pprintTitle{%
  \let\@oddhead\@empty
  \let\@evenhead\@empty
  \def\@oddfoot{\reset@font\hfil\thepage\hfil}
  \let\@evenfoot\@oddfoot
}
\makeatother

\allowdisplaybreaks[1] 

\begin{document} 


\begin{frontmatter}


\title{Study of a fifth order PDE using symmetries}
\author{Stylianos Dimas and Igor Leite Freire}
\address{\it Centro de Matem\'atica, Computa\c c\~ao e Cogni\c c\~ao\\ \it Universidade Federal do ABC - UFABC\\ \it
Avenida dos Estados, 5001, Bairro Bangu
$09.210-580$\\\it Santo Andr\'e, SP - Brasil\\
\rm E-mails: s.dimas@ufabc.edu.br and igor.freire@ufabc.edu.br/igor.leite.freire@gmail.com
}
\begin{abstract}
We study a family of PDEs, which was derived as an approximation of an extended Lotka--Volterra system, from the point of view of symmetries. Also, by performing the self adjoint classification on that family we offer special cases possessing non trivial conservation laws. Using both classifications we justify the particular cases studied in the literature, and we give additional cases that may be of importance. 
\end{abstract}
\end{frontmatter}

\section{Introduction}

In a recent paper, Zilburg and Rosenau considered the  equation
\bb\label{1.1}
	u_t=p_1(u^2)_x + h^2 p_3 (u u_{xx}+\beta u_x^2)_x+h^4\left[q_0uu_{5x}+q_1u_xu_{4x}+q_2u_{2x}u_{3x}\right]
\ee
derived by taking an approximation of an extended Lotka--Voltera (ELV) system up to the $\mathcal{O}(h^6)$ order, \cite{ZiRo2k16}. Equation (\ref{1.1}) contains, among other models, the celebrated $K(2,2)$ equation\footnote{$q_0=q_1=q_2=0$, $\beta=1$ and $h^2p_3=2$} --- also known as Rosenau--Hyman equation --- which is best known for admitting compacton solutions, \cite{RosHy93}. Moreover, equation (\ref{1.1}) may be capable to describing higher order nonlinear excitations observed in numerical simulations that equations like the KdV fail to replicate, see \cite{Ro2k16} and referenced therein. 

Apart from the Rosenau--Hyman equation, another important special case comes up if we take the choice $q_0=q_1=\beta=0$. Now equation (\ref{1.1}) can be rewritten as
\bb\label{1.3}
u_t=\f{\p}{\p x}\left[bu^2+2\kappa uu_{xx}+(u_{xx})^2\right],
\ee     
where $t\mapsto 2ht/q_2$, $x\mapsto xh, b=2p_1/q_2$ and $\kappa=p_3/q_2$. And if 
\bb\label{1.4}
b=\kappa^2,
\ee
equation (\ref{1.3}) becomes
\bb\label{1.5}
u_t=\f{\p}{\p x}\left[{\cal L}u\right]^2,
\ee 
where ${\cal L}$ is the Schr\"ondinger operator with constant potential $\kappa$,
\begin{equation*}
{\cal L}=\f{\p^2}{\p x^2}+\kappa.
\end{equation*}

In this paper we shall cast some light into structural properties of equation (\ref{1.1}) and improve our knowledge on it, complementing the results obtained in \cite{ZiRo2k16}. Along with this desire, we shall investigate the invariance properties of equation \eqref{1.1}. First we carry out its Lie group classification and we show how all these special cases found in the literature occur naturally from it. In addition, we give some models of potential interest not found in the literature along with group invariant solutions. This is done in section \ref{group} using the symbolic package SYM for Mathematica\texttrademark{} developed by SD \cite{DiTs2k5a}.

Once we have the group classification, and hence the Lie Algebra for each possible case, in section \ref{laws}, we turn our attention to the self adjoint classification of equation \eqref{1.1}. Applying the techniques developed in \cite{Ibra2k7a, Ibr2k11a}, again with the invaluable help of SYM, we unearth special cases possessing non trivial conservation laws. Next, in section \ref{remarks}, we discuss some properties of the solutions of (\ref{1.1}) from the conservation laws established in section \ref{laws}. Our conclusions are presented in section \ref{conclusion}.

%
%
%
%
%
%
%
%
%
%

\section{Point symmetries and group invariant solutions}\label{group}
First of all, utilizing the equivalent transformation machinery \cite{Ibra2k9}, we can eliminate  the parameter $h$ by using the transformation $(x,t,u)\rightarrow(\frac{x}{h},t,\frac{u}{h})$. Hence, the comment made by the authors in \cite{ZiRo2k16} that the equation~\eqref{1.1} is quasi continuum --- at least from the mathematical standpoint --- is not accurate: up to the cutoff at $\mathcal{O}(h^4)$ the equation do not carry a trace of its discrete origins. Therefore, from this point on we set that $h=1$. We proceed now with the Lie group classification of equation~\eqref{1.1}. 

\subsection{Group classification}

By using SYM interactively we arrive to the group classification summarized in table~\ref{tab1}:

\begin{table}[htp]
\begin{center}\scriptsize
\begin{tabular}{c|cccccc||lc}\hline
Case & $p_1$ & $p_3$ & $\beta$ & $q_0$ & $q_1$ & $q_2$ & Form & Symmetries ($\Span$)\\ \hline
0 & \multicolumn{5}{c}{general case} & &  \eqref{1.1}      & $\mathfrak{X}_1,\,\mathfrak{X}_2,\,\mathfrak{X}_3$\\
1 & $0$ & $0$ & $\forall$ & $\ne0$ & $\forall$ & $\forall$ &  \begin{tabular}{@{}l@{}} $u_t=q_0uu_{5x}+q_1u_xu_{4x}+$\\ $\qquad\quad q_2u_{2x}u_{3x}$ \end{tabular} & $\mathfrak{X}_4=x\p_x+5u\p_u$ \\[10pt]
2 & $\ne0$ & $0$ & $\forall$ & $0$ & $\ne0$ & $\forall$ & \begin{tabular}{@{}l@{}} $u_t=p_1(u^2)_x +q_1u_xu_{4x}+$\\ $\qquad\quad q_2u_{2x}u_{3x}$ \end{tabular} & $\mathfrak{X}_4=2p_1t\p_x-\p_u$ \\[10pt] 
3 & $0$ & $\ne0$ & $\forall$ & $0$ & $\ne0$ & $0$ & $u_t= p_3 (u u_{xx}+\beta u_x^2)_x$ & $\mathfrak{X}_4=x\p_x+3u\p_u$ \\[5pt]
4 & $0$ & $0$ & $\forall$ & $0$ & $\ne0$ & $\forall$ & $u_t=q_1u_xu_{4x}+q_2u_{2x}u_{3x}$ & \begin{tabular}{@{}l@{}} $\mathfrak{X}_4=\p_u,$\\ $\mathfrak{X}_5=x\p_x+5t\p_t$ \end{tabular}\\[10pt]
5a & $=\frac{p_3^2}{2q_2}$ & $\ne0$ & $0$ & $0$ & $0$ & $\ne0$ & \begin{tabular}{@{}l@{}} $u_t=\frac{p_3^2}{2q_2}(u^2)_x +  p_3 (u u_{xx})_x+$\\ $\qquad\quad q_2u_{2x}u_{3x},\, \kappa>0$\end{tabular} & \begin{tabular}{@{}l@{}} $\mathfrak{X}_4=\cos\left(\sqrt{\kappa}\,x\right)\partial_u,$\\ $\mathfrak{X}_5=\sin\left(\sqrt{\kappa}\,x\right)\partial_u$ \end{tabular} \\[10pt]
5b & $=\frac{p_3^2}{2q_2}$ & $\ne0$ & $0$ & $0$ & $0$ & $\ne0$ & \begin{tabular}{@{}l@{}} $u_t=\frac{p_3^2}{2q_2}(u^2)_x +  p_3 (u u_{xx})_x+$\\ $\qquad\quad q_2u_{2x}u_{3x},\, \kappa<0$\end{tabular} & \begin{tabular}{@{}l@{}} $\mathfrak{X}_4=e^{\sqrt{\lvert\kappa\rvert}x}\partial_u,$\\ $\mathfrak{X}_5=e^{-\sqrt{\lvert\kappa\rvert}x}\partial_u$ \end{tabular} \\[10pt]
6 & $\forall$ & $0$ & $\forall$ & $0$ & $0$ & $0$ & $u_t=p_1(u^2)_x$      & $\infty$\tablefootnote{In this case we have  a quasi linear PDE that can be easily solved analytically, so there is no need to explicitly present its symmetry algebra here.}\\ 
\end{tabular}
\end{center}
\caption{Group Classification of \eqref{1.1} where $\mathfrak{X}_1=\partial_x,\,\,\mathfrak{X}_2=\partial_t,\,\,\mathfrak{X}_3=t\partial_t-u\partial_u$ and $\kappa=p_3/q_2$.}\label{tab1}
\end{table}%

It is evident that the case studied \emph{ad hoc} in \cite{ZiRo2k16}, namely the cases $5a$ and $5b$ in our classification, are featured by the group classification admitting an exceptional high dimensional Lie algebra in relation with the general case.  Besides that, case $4$ also admits a $5$ dimensional Lie algebra as in the cases $5a$ and $5b$. A fact that suggests that it might be of physical interest too. We continue by providing invariant solutions for these three cases.

\subsection{Invariant solutions for case $5$}\label{sol1}

The case $5$ corresponds to the family of PDEs\footnote{after applying the equivalence transformation $(x,t,u)\rightarrow(\sqrt{\lvert\kappa \rvert}x,t,q_2\sqrt{\lvert\kappa\rvert^5} u)$.}
\bb\label{2.2.1}
u_t=\pm\frac{1}{2}(u^2)_x +  (u u_{xx})_x+u_{2x}u_{3x}
\ee
that admits a $5$--dimensional Lie algebra which is represented by the span of two different sets of five infinitesimal generators depending on the sign of the term $(u^2)_x$, see table~\ref{tab1}. Since we are dealing with the same Lie algebra the algebraic part of the analysis that follows is the same and the difference occurs only when we write down the similarity solutions. 

As we intend to give the invariant solutions for equation~\eqref{2.2.1} it is advantageous first to construct the optimal set of the corresponding $5$--dimensional Lie algebra, see \cite{Hydon2k, Olver2k} for more details.   That is: 
\begin{align*}
	&\mathfrak{X}_1&&\mathfrak{X}_2&&\mathfrak{X}_3&&\mathfrak{X}_4&&\mathfrak{X}_5\\
	&\mathfrak{X}_1+\mathfrak{X}_2&&\mathfrak{X}_1+\lambda\mathfrak{X}_3&&\mathfrak{X}_4+\mathfrak{X}_2&&
	&&
\end{align*}

The invariance under $x$ translations ($\mathfrak{X}_1$) leads to constant solutions, while the invariance under translations in $t$ ($\mathfrak{X}_2$) may either provide exponential or sinusoidal solutions of $x$, depending on the sign. For the remaining symmetries of the optimal set we get the following similarity solutions:
{\renewcommand{\thefootnote}{\fnsymbol{footnote}}
\begin{align*}
	&                   && +1 (\kappa>0) && -1 (\kappa <0)\\
	&\mathfrak{X}_3:&& u= \frac{\phi(x)}{t}&& u= \frac{\phi(x)}{t}\\
	&\mathfrak{X}_1+\mathfrak{X}_2:&& u= \phi(x-t)\text{{}\footnotemark[2]}
	&& u= \phi(x-t)\text{{}\footnotemark[2]}\\
	&\mathfrak{X}_1+\lambda\mathfrak{X}_3:&& u= e^{-\lambda x}\phi(e^{-\lambda x}t)&& u= e^{-\lambda x}\phi(e^{-\lambda x}t)\\
	&\mathfrak{X}_4+\mathfrak{X}_2:&& u= \cos\left(x\right)t+\phi(x)&& u=  e^{x} t+\phi(x)\\
\end{align*} \footnotetext[2]{note the discrete symmetry $(t,u)\rightarrow(-t,-u)$}
}
We note at this point that:
\begin{itemize}
	\item for the symmetry $\mathfrak{X}_1+\mathfrak{X}_2$ we get the type of solutions --- modulo an equivalence transformation --- that Rosenau et. al. studied in \cite{ZiRo2k16}. 
	\item To obtain the exact form of the function $\phi$ we need to substitute each similarity solution to equation~\eqref{2.2.1} and resolve the reduced equation that emerge. 
\end{itemize}

\subsection{Invariant solutions for case $4$}\label{sol2}

Now we turn our attention to a special case of equation~\eqref{1.1} that is not mentioned in \cite{ZiRo2k16}, \emph{videlicet}
\begin{equation}\label{eq:case4}
	u_t=u_xu_{4x}+qu_{2x}u_{3x},
\end{equation}
where $q=\frac{q_2}{q_1}$\footnote{after applying the equivalence transformation $t\rightarrow q_1 t$}. For this case, the optimal set of the Lie algebra admitted by equation~\eqref{eq:case4} is :
\begin{align*}
	&\mathfrak{X}_1&&\mathfrak{X}_2&&\mathfrak{X}_4&&\mathfrak{X}_5\\
	&\mathfrak{X}_3+\lambda\mathfrak{X}_5&&\mathfrak{X}_1+ \mathfrak{X}_3&&\mathfrak{X}_4+\mathfrak{X}_5&&\mathfrak{X}_2+\mathfrak{X}_4&&\\
	&\mathfrak{X}_1+\mathfrak{X}_2&&\mathfrak{X}_2+\mathfrak{X}_3-\frac{1}{5}\mathfrak{X}_5 
\end{align*}

Similarly to the previous cases, the invariance under translations  in $x$  and $t$  lead to constant solutions and  implicit solutions involving Hypergeometric functions respectively. The third symmetry of the optimal set, $\mathfrak{X}_4$, only denotes that the function itself do not appear in \eqref{eq:case4}  and as such gives no invariant solutions\footnote{It affirms the fact that if $u$ is a solution so is $u+c$ where $c$ a constant.}. For the remaining symmetries of the optimal set we have the following similarity solutions:
\begin{align*}
	\mathfrak{X}_5:& u= \phi\left(\frac{x^5}{t}\right)\\
	\mathfrak{X}_3+\lambda \mathfrak{X}_5:& u= \left\{\begin{aligned} &x^{-\frac{1}{\lambda}}\phi\left(\frac{x^{5+\frac{1}{\lambda}}}{t}\right),& \lambda\ne0\\ &\frac{\phi(x)}{t},& \lambda=0\end{aligned}\right.\\
	\mathfrak{X}_1+\mathfrak{X}_3:& u= e^{-x}\phi\left(e^{-x}t\right)\\
	\mathfrak{X}_4+ \mathfrak{X}_5:& u= \phi\left(\frac{x^5}{t}\right)+\ln x\\
	\mathfrak{X}_2+\lambda \mathfrak{X}_4:& u= \phi(x)+t\\
	\mathfrak{X}_1+\lambda \mathfrak{X}_2:& u= \phi(x-t)\\
	\mathfrak{X}_2+\mathfrak{X}_3-\frac{1}{5}\mathfrak{X}_5:& u= x^5\phi\left(x^5e^t\right)
\end{align*}

\section{Self adjoint classification and conservation laws derived from point symmetries}\label{laws}


In addition to the group classification we performed the self adjoint classification of the PDE, that is, to reveal the special cases that are (strictly, quasi and nonlinearly) self adjoint. Ibragimov proposed in \cite{Ibra2k7a} a procedure that utilize that property to construct conservation laws by introducing the concept of a \emph{formal Lagrangian} $\mathcal{L}$. A formal Lagrangian is nothing more than the differential equation multiplied by a new dependent variable $v(x,y)$, \emph{viz.}
\begin{multline*}
\mathcal{L}=v \Delta = \\
					v (u_t-p_1(u^2)_x-h^2 p_3 (u u_{xx}+\beta u_x^2)_x-h^4\left[q_0uu_{5x}+q_1u_xu_{4x}+q_2u_{2x}u_{3x}\right]).
\end{multline*}
When 
$$
\frac{\delta}{\delta u}(\mathcal L)= \lambda\Delta,
$$
where $\frac{\delta}{\delta u}$ the Euler operator, our equation is self adjoint. It is clear that the self adjointness depends on the choice of the dependent variable $v$:
\begin{itemize}
	\item if $v=u$ the equation is \emph{strictly} self adjoint;
	\item if $v=\Phi(u)$ the equation is \emph{quasi} self adjoint;
	\item if $v=\Phi(x,t,u)$ the equation is \emph{nonlinearly} self adjoint, and finally
	\item if $v=\Phi(x,t,u,u_x,u_{xx},\dots)$ the equation is \emph{generalized nonlinearly} self adjoint.
\end{itemize}
The fact that a differential equation is self adjoint means that its symmetries are also variational symmetries. Therefore, we can use the Noether theorem, expressed as an operator identity,
$$
	\mathfrak{X} + D_i(\xi^i)=W^\alpha\frac{\delta}{\delta u^\alpha}+D_i\mathcal{N}^i,
$$   
where $\mathfrak{X}=\xi^i\frac{\partial}{\partial x^i}+\eta^\alpha\frac{\partial}{\partial u^\alpha}$ the chosen symmetry appropriately prolonged and $W^\alpha=\eta^\alpha- \xi^ju^\alpha_j$, to obtain the corresponding conserved vector $\mathcal N(\mathcal L)$. 

By solving the identity for $\mathcal N$ we obtain the formula:
\begin{equation*}\label{eq:ConservationVectorFormula}
	\mathcal{N}^i\mathcal{L}=\xi^i\mathcal L + \sum_{i_1+\cdots+i_n=0}^\infty D_{i_1}\cdots D_{i_n}(W^\alpha) \frac{\delta^*\mathcal L}{\delta^* u^\alpha_{i\ i_1\ \dots\ i_n}},
\end{equation*}
with  
\begin{equation*}
	\frac{\delta^*}{\delta^* u^\alpha_{i_1\dots i_n}}=\frac{\partial}{\partial u^\alpha}+\sum_{s=j_1+\cdots+j_n=1}^\infty(-1)^{s} \frac{\binom{s}{j_1,\dots,j_n}}{\binom{s+i_1+\cdots i_n}{i_1+j_1,\dots,i_n+j_n}}D_{j_1}\cdots D_{j_n}\frac{\partial}{\partial u^\alpha_{(i_1+j_1)\,\dots\,(i_n+j_n)}},
\end{equation*}
where $\binom{N}{i_1,i_2,\dots,i_r}=\frac{N}{i_1!i_2!\dots i_r!},\, N=i_1+i_2+\cdots+i_r$ is the multinomial and $i_j\ge0$ denotes the order of the derivative for the $j^{th}$ independent variable. For a more detailed survey on the aforementioned concepts, fully integrated in SYM, see also \cite{Ibra2k7a, Ibr2k11a}.

\subsection{The strictly self adjoint cases}
The cases $q_2=5(q_1-2q_0), p_3=0$ and $q_2=5(q_1-2q_0), \beta=\frac{1}{2}$, namely the PDEs
\begin{align}
	&u_t=&&p_1(u^2)_x +\left[q_0uu_{5x}+q_1u_xu_{4x}+5(q_1-2q_0)u_{2x}u_{3x}\right]=\notag\\
	&&&p_1(u^2)_x +\left(q_0u u_{4x}+(q_1-q0)u_xu_{3x}+(2q_1-\frac{9}{2}q_0)u_{xx}^2\right)_x\label{eq:strict1}\\ 
	\intertext{and}
	&u_t=&&p_1(u^2)_x +  p_3 (u u_{xx}+\frac{1}{2}u_x^2)_x+\left[q_0uu_{5x}+q_1u_xu_{4x}+5(q_1-2q_0)u_{2x}u_{3x}\right]=\notag\\
	&&&p_1(u^2)_x +  p_3 (u u_{xx}+\frac{1}{2}u_x^2)_x+ \\
	&&&\qquad\qquad\qquad\left(q_0u u_{4x}+(q_1-q0)u_xu_{3x}+(2q_1-\frac{9}{2}q_0)u_{xx}^2\right)_x,\notag\label{eq:strict2}
\end{align}
are strictly self adjoint. Observe that both cases can be written as conservation laws, the reason will be revealed, and proved,  in the following section.
 
\subsection{The quasi self adjoint cases}

We turn now to the quasi self adjoint case. 
 The result of this kind of classification is gathered in table~\ref{selfAdjointClass}.

\begin{table}[htp]
\begin{center}
\begin{tabular}{c|cccccc||lc}\hline\label{tab2}
Case & $p_1$ & $p_3$ & $\beta$ & $q_0$ & $q_1$ & $q_2$ & $\Phi$\\ \hline
1 & $\forall$ & $0$ & $\forall$ & $0$ & $0$ & $0$ & $\Phi(u)$  \\
2 & $\forall$ & $0$ & $\forall$ & $\ne0$ & $\ne q_0$ & $5(q_1-q_0)$ & $c_1+c_2\ln u$  \\
3 & $\forall$ & $\ne0$ & $0$ & $\forall$ & $q_0$ & $0$ & $c_1+c_2\ln u$  \\
4 & $\forall$ & $0$ & $\forall$ & $\ne0$ & $\ne q_0$ & $0$ & $c_1+c_2u^{\frac{q_1}{q_0}-1}$  \\
5 & $\forall$ & $\ne0$ & $\ne0$ & $\forall$ & $(1+2\beta)q_0$ & $0$ & $c_1+c_2u^{2\beta}$  \\
6 & $\forall$ & $0$ & $\forall$ & $0$ & $\ne0$ & $5q_1$ & $c_1+c_2u+c_3u^2$  \\
7 & $\forall$ & $0$ & $\forall$ & $\ne0$ & $\forall$ & $5(q_1-3q_0)$ & $c_1+c_2u^2$  \\
8 & $\forall$ & $\ne0$ & $1$ & $\forall$ & $\forall$ & $5(q_1-3q_0)$ & $c_1+c_2u^2$  \\
9 & $\forall$ & $0$ & $\forall$ & $0$ & $\forall$ & $\ne5q_1$ & $c_1$  \\
10 & $\forall$ & $0$ & $\forall$ & $\ne0$ & $q_0$ & $0$ & $c_1$  \\
11 & $\forall$ & $\ne$ & $0$ & $\forall$ & $\forall$ & $\ne0$ & $c_1$  \\[10pt]
12 & $\forall$ & $0$ & $\forall$ & $\ne0$ & $\forall$ &\begin{tabular}{@{}l@{}}  $\ne0$\\and \\$\ne5(q_1-q_0)$\\ and\\ $\ne5(q_1-2q_0)$\\and\\ $\ne5(q_1-3q_0)$ \end{tabular} & $c_1$  \\[50pt]
13 & $\forall$ & $0$ & $\forall$ & $\ne0$ & $\forall$ & $5(q_1-2q_0)$ & $c_1+c_2 u$  \\
14 & $\forall$ & $\ne0$ & $\frac{1}{2}$ & $\forall$ & $\forall$ & $5(q_1-2q_0)$ & $c_1+c_2u$  \\
\end{tabular}
\end{center}
\caption{Quasi self adjoint classification of equation \eqref{1.1}, $\Phi(u)$ is an arbitrary function of $u$}\label{selfAdjointClass}
\end{table}%

We observe that:
\begin{itemize}
	\item In any case $v$ can be a constant. This means that the equation~\eqref{1.1} can be written in the form $u_t=(\cdot)_x$. Verily,
		\begin{align*}
			u_t&=p_1(u^2)_x+ \left( h^2u(h^2 q_0 u_{xx}+ p_3 u)_{xx}\right)_x+ \frac{(q_0-q_1+q_2)}{2}h^4\left( u_{xx}^2\right)_x+\\
   				&\quad\  h^2 \left(u_{x} (\beta  p_3 u-h^2 (q_0-q_1)  u_{xx})_x\right)_x.
		\end{align*}
	\item The special case briefly discussed in \cite{ZiRo2k16},
	$$
		u_t=p_1(u^2)_x + p_3 (u u_{xx})_x+q_0\left(uu_{4x}\right)_x,
	$$
	corresponds to our case 3 in table~\ref{selfAdjointClass}. And apart from the obvious conservation law --- using either the  symmetry $\mathfrak{X}_2$ or $\mathfrak{X}_2$ and $\Phi=1$ --- it has also another one. Namely,
	\begin{multline*}
		(\log (u) u)_t-\left(\frac{1}{2} \left(2  (\log (u)+1) u\left(p_3 u_{xx}+q_0 u_{xxxx}\right)+ \right.\right.\\
			\left.p_1  (2 \log
  		 	(u)+1)u^2- u_{x} \left(p_3 u_{x}+2 q_0 u_{xxx}\right)+q_0 u_{xx}{}^2\right)\biggr)_x=0.
	\end{multline*}
	obtained with the symmetry $\mathfrak{X}_3$ and $\Phi=\log u$.
	\item The last two cases correspond to the strictly self adjoint ones and the existence of the constant explains the way we have written them in the previous section. 
\end{itemize}

\subsection{The generalized nonlinearly self adjoint cases}

Beyond the special cases we illustrated so far 
equation~\eqref{1.1}  possesses a case that has the rare property to be generalized nonlinearly  self adjoint. That happens  when $p_1=0,\beta=-\frac{1}{4}\text{ and }q_1=\frac{3}{2}q_0\ne0$, that is
$$
u_t=
	†\frac{1}{4} \left(4 u \left(p_3 u_{xx}+q_0 u_{xxxx}\right)-p_3 u_{x}{}^2-(q_0-2 q_2)u_{xx}{}^2+2 q_0 u_{xxx} u_{x}\right)_x.
$$
This case admits the substitution $\Phi=u_{xx}$ which in combination with the symmetry $\mathfrak{X}_3$ yields the non trivial conservation law:
\begin{multline*}
	\left(u_x^2\right)_t+\left(u\left(p_3u_{xx}^2-q_0u_{xxx}^2+2q_0u_{xx}u_{xxxx}\right) +\right.\\ 
		q_0u_xu_{xx}u_{xxx}+\frac{2q_2-q_0}{3}u_{xx}^3-2u_xu_t\biggr)_x=0.
\end{multline*}


\section{Remarks on the constants of motion}\label{remarks}

Here we present some facts regarding the conservation laws established in the previous section.
\begin{enumerate}
\item Since equation (\ref{1.1}) is itself a conservation law, this implies that the quantity
\bb\label{4.1}
{\cal H}_0[u]=\int^{+\infty}_{-\infty}u\,dx
\ee
is a constant of motion. Particularly, if $u$ is a non-negative function, (\ref{4.1}) implies on the conservation of the $L^1(\R)-$norm of the solutions of (\ref{1.1}) rapidly decaying to $0$, jointly with its derivatives, at the infinity.

\item In the strictly self adjoint cases (cases 6, 13 and 14 in Table~\ref{selfAdjointClass}) we have the conserved quantity $u^2$, which implies that
\bb\label{4.2}
{\cal H}_1[u]=\int^{+\infty}_{-\infty}u^2\,dx
\ee
is also a constant of motion, which mathematically corresponds to the existence of square integrable solutions of (\ref{1.1}). 

\item Cases 2 and 3 in Table~\ref{selfAdjointClass} have a constant of motion of the logarithmic type:
\bb\label{4.3}
{\cal H}_2[u]=\int^{+\infty}_{-\infty}u\,\ln{|u|}\,dx.
\ee


\item For the case $p_1=0,\beta=-\frac{1}{4}\text{ and }q_1=\frac{3}{2}q_0\ne0$, in addition to (\ref{4.1}), we have the constant of motion
\bb\label{4.4}
\cal{H}_3=\int^{+\infty}_{-\infty} u_x^2\,dx.
\ee

In addition, if $q_2=-5q_0/2$ the PDE is also quasi self adjoint, falling into the case 13, hence admitting also  (\ref{4.2}) as a constant of motion on the (rapidly decaying) solutions of (\ref{1.1}). Combining (\ref{4.2}) and (\ref{4.4}), we obtain the conservation of the $H^1(\R)-$norm of the solutions of (\ref{1.1}) satisfying these constraints.

\item Cases 4, 5, 6, 7 and 8 have the constant of motion
\bb\label{4.5}
{\cal H}_4[u]=\int^{+\infty}_{-\infty}u^{\sigma}dx,\ \sigma\ne0
\ee
where $\sigma=q_1/q_2$ (case 4), $\sigma=2\beta+1$ (case 5) and $\sigma=3$ for the remaining cases. 
\end{enumerate}

\section{Conclusion}\label{conclusion}

In the present work we showed how symmetries can help in a systematic and thorough study of a family of nonlinear PDEs of the fifth order. By its group classification, and then, its self adjoint classification we were able not only to retrieve the cases that Zilburg and Rosenau studied \emph{ad hoc} but also to give  additional cases possessing non trivial, and not obvious, conservation laws. Also we proved that equation~\eqref{1.1} can be written in the form $u_t=(\cdot)_x$ for any value of the parameters. A bit of information that greatly helps the study of equation~\eqref{1.1} that Zilburg and Rosenau performed for only two special cases in \cite{ZiRo2k16}  commenting that ``in its full glory appears to be well beyond our ability to analyze it''. Actually, this very comment was our main motivation for studying this class of PDEs by looking into their structure --- their \emph{DNA} in a matter of speaking --- its symmetries. From them, we were able to construct solutions and conservation laws of (\ref{1.1}), another fact that shows the usefulness of the symmetry analysis we performed. Furthermore, we uncovered a very interested case, namely the case  $p_1=0,\beta=-\frac{1}{4}\text{ and }q_1=\frac{3}{2}q_0\ne0$, where we have a generalized nonlinear  self adjoint PDE which yields a conservation law with a higher order characteristic.

In a future work we will utilize the symmetry machinery in order to study systematically higher order expansion cutoffs of the ELV system.

\section*{Acknowledgements} 
The work of I. L. Freire is partially supported by FAPESP (grant no. 2014/05024-8) and CNPq (grant no. 308941/2013-6).

\section*{References}

\bibliographystyle{plain}
\bibliography{Bibliography}

\end{document}